\documentclass[a4paper,11pt]{article}
\usepackage{fullpage}
\usepackage{graphicx}
\usepackage{amsmath,amsfonts,mathrsfs,amssymb,latexsym}
\usepackage{theorem}
\usepackage{enumerate}
\usepackage[dvips]{hyperref}

\def\Q{\ensuremath{\mathbb{Q}}}
\def\N{\ensuremath{\mathbb{N}}}

\def\F{\ensuremath{\mathbb{F}}}

\def\H{\ensuremath{\mathbb{H}}}
\def\Z{\ensuremath{\mathbb{Z}}}

\def\Pp{\ensuremath{\mathscr{P}}}

\def\Y{\ensuremath{\mathscr{Y}}}
\def\X{\ensuremath{\mathscr{X}}}

\def\l{\ensuremath{\langle}}
\def\r{\ensuremath{\rangle}}

\def\Dd{\ensuremath{\mathscr{D}}}

\def\Zz{\ensuremath{\mathcal{Z}}}

\newcommand{\Jacobi}[2]{\ensuremath{\left(\frac{#1}{#2}\right)}}

\newcommand{\ii}{\ensuremath{\mathsf{i}}}
\newcommand{\jj}{\ensuremath{\mathsf{j}}}
\newcommand{\kk}{\ensuremath{\mathsf{k}}}

\newcommand{\girth}{\ensuremath{\mathrm{girth}}}
\newcommand{\conj}[1]{\ensuremath{\overline{#1}}}
\newcommand{\Cayl}{\ensuremath{\mathscr{C}\! ay}}

\def\proof{\noindent {\sc Proof:~}}
\def\foorp{\hfill {$\blacksquare$}\medskip}

\newtheorem{Theo}{Theorem}[section]
\newtheorem{Lem}[Theo]{Lemma}

\newtheorem{Prop}[Theo]{Proposition}
\newtheorem{Def}[Theo]{Definition}
{
\theorembodyfont{\upshape}
\newtheorem{Rem}[Theo]{Remark}
}
\title{Regular graphs of large girth and arbitrary degree}
\author{X.~Dahan\footnote{Supported by the GCOE Project ``Math-for-Industry''
of Kyushu university}\\
ISEE, Department of Informatics, Kyushu university, Japan\\
{\tt xdahan at gmail.com}
}
\date{}
\begin{document}
\maketitle

\begin{abstract}

For every integer $d \ge 10$, we construct infinite families 
$\{G_n\}_{n \in \N}$  of $d + 1$-regular graphs
which ave a large girth $\ge \log_d
|G_n|$, and for $d$ large enough $\ge 1, 33 \cdot \log_d
|G_n|$. These are
Cayley graphs on $PGL_2(\F_q)$ for a special set of $d + 1$ 
generators whose choice is related to
the arithmetic of integral quaternions. 
These graphs are inspired by the Ramanujan graphs
of Lubotzky-Philips-Sarnak and Margulis,
 with which they coincide when $d$ is prime. 
When $d$ is {\em not} equal to the power of an odd prime, 
this improves the previous construction of
Imrich in 1984 where he obtained infinite families $\{I_n\}_{n\in \N}$
of $d + 1$-regular graphs, realized
as Cayley graphs on $SL_2(\F_q)$, and which are displaying a girth 
$\ge 0, 48 \cdot \log_d |I_n|$. And
when $d$ is equal to a power of 2, 
this improves a construction by Morgenstern in 1994 where
certain families $\{M_n\}_{n \in N}$
of $2^k+1$-regular graphs were shown to have girth 
$ \ge  2/3 \cdot \log_{2^k} |M_n|$.
\end{abstract}

\section{Introduction}

The ``Moore bound'' follows from a simple counting argument,
and permits to show that a $d$-regular graph $G$ of order $|G|$,
admits the following upper bound on its girth
(see~\cite[Ch.~III, Theorem~1.2]{Bo78}:
\begin{equation}
\label{eq:Moore1}
\girth(G) \le  \begin{cases}
2 \log_{d-1} |G| + 1& \text{if $\girth(G)$ is {\em odd},}\\
2 \log_{d-1} |G| + 2 - 2\log_{d-1} 2 & \text{if $\girth(G)$ is {\em even}.}
\end{cases}
\end{equation}
This implies that for $d\ge 5$, 
\begin{equation}\label{eq:Moore2}
\girth(G) \le ( 2+ \frac 2 {\log_{d-1} |G|}) \log_{d-1} |G|. 
\end{equation}
It is not known if this bound is tight.
A convenient way to formulate what is meant by ``tight'', is to consider
large graphs, and even better, infinite family of constant degree
regular graphs.
Let us recall the following definition: a family of $d$-regular graphs $\{G_n\}_{n\in \N}$
is said to have {\em large girth} if there exists
a constant $c>0$ independent of $n$ (but possibly dependent on $d$), such that:
\begin{equation*}
\girth(G_n) \ge (c + o_n(1)) \log_{d-1} |G_n|.
\end{equation*}

The property of large girth, besides its
own theoretical interest, can be applied to
 LDPC codes. This approach was pioneered by Margulis in~\cite{Ma82}, where he gave the first constructive example of
 a family of LDPC codes of unbounded minimum distance
 by providing explicit families of regular graphs of large girth. Another important application of large girth
graphs to LDPC codes can be found in~\cite{Ta81}.

Given an infinite family of $d$-regular graph, let us define:
$$
\gamma (\{G_n\}) = \liminf_{n\rightarrow \infty}
\frac{\girth(G_n)}{\log_{d-1} |G_n|},
\quad \text{and}\quad 
\gamma_d:=\sup_{\{G_n\} \text{ family of $d$-regular graphs}}
\gamma (\{G_n\}).
$$
What the bound~\eqref{eq:Moore2} says is that $\gamma_d\le 2$, for any $d\ge 3$.
As for lower bound, it was proved that $\gamma_d\ge 1$
by Erd\"os and Sachs~\cite{ErSa63} for any $d\ge 3$.
Their proof, of probabilistic nature, did not
provide explicit families $\{G_n\}_n$.
Currently, the best lower bounds for $\gamma_{d}$
that are deduced from {\em explicit} examples of family of graphs,
are:
\begin{enumerate}
\item $\gamma_{d} \ge \frac  4 3$ for $d=p^k+1$, $p$  an odd prime
and $k\in\N^\star$, (for $d=p+1$
where $p$ is an
odd prime; this was first achieved by Lubotzky-Philips-Sarnak~\cite{LPS88}
and independently by Margulis~\cite{Ma88},
then later also by Lazebnik-Ustimenko~\cite{LaUs95} with a different construction.
Finally, Morgenstern~\cite{Mo94} treated the case
$d-1$ equal to any prime power).
\item $\gamma_d \ge \frac 2 3$ 
  for $d = 2^k + 1$ with $k\in \N^\star$.
This is also due to Morgenstern \cite[Theorem 5.13-3]{DaSaVa03}.
\item $\gamma_{d} \ge 0,48$ for other values of $d$ (this is due to
  Imrich~\cite{Im84}, extending the method of 
Margulis~\cite{Ma82} 
where it was proved that $\gamma_{d}\ge \frac 4 9$
for odd $d$).
\end{enumerate}
These are the best results we are aware of.
This paper presents improvements on the lower bounds
on $\gamma_{d}$ in the cases 2 and 3, that is, when $d-1$ is not a prime power.
For other values of $d$, the lower bounds that would be obtained
do not improve the best ones shown in the case 1.
That is why we focus only on the cases where 
$d-1$ is not the power of an odd prime,
and henceforth consider only $d\ge 10$ (lower values are either
prime powers or non manageable by our method). 
\begin{Theo}\label{th:contrib}
For any integer $d\ge 10$,
which is not a prime power, there is an explicit infinite family
$\{G_n\}_n$ of $d+1$-regular graphs, bipartite and connected,
as well as having large girth. Precisely:
\begin{equation}
\girth(G_n) \ge c(d) \log_d |G_n| - \log_d 4,
\end{equation}
where $c(d)$ is a constant independent of $n$,
such that $c(d)\le \frac 4 3$ and:
\renewcommand{\arraystretch}{1.5}
$$\text{case $d$ odd} \left\{\begin{array}{lccl}
  \label{eq:1,33}\text{if $d\ge 1335$,} & c(d) & \ge & 1,33 \\
 \label{eq:1,3}\text{if $35 \le d \le 1331$} & c(d) & \ge & 1,3 \\
 \label{eq:1,24}\text{if $15\le d \le 31$,} & c(d) & \ge & 1,27 
\end{array}\right.
$$
$$\hspace{-6cm}\text{case $d$ even} \left\{\begin{array}{lccl}
  \label{eq2:1,33}\text{if $d\ge 4826$,} & c(d) & \ge & 1,33 \\
 \label{eq2:1,3}\text{if $184 \le d \le 4824$} & c(d) & \ge & 1,3 \\
 \label{eq2:1,25}\text{if $44 \le d \le 182$,} & c(d) & \ge & 1,25\\
  \label{eq2:1,1}\text{if $22 \le d \le 42$,} & c(d) & \ge & 1,1\\
\label{eq2:misc}\lefteqn{c(10)\ge 1,28\quad c(12)\ge 1,12 \quad c(14)\ge 1,19
\quad c(18)\ge 1,3\quad  c(20)\ge 1,061.} 
\end{array}\right.
$$
\medskip

Related to the families $\{G_n\}_n$,
there are also explicit families of $d+1$-regular
graphs $\{H_n\}_n$, connected and {\em non-bipartite}, for which the
girth verifies:
\begin{equation*}\label{eq3:main}
\girth( H_n ) \ge \frac {c(d)} 2 \log_{d} |H_n|.
\end{equation*}
\end{Theo}
The family $\{G_n\}_n$ will be $\X_d$ and $\{H_n\}_n$
will be $\Y_d$ introduced in Definition~\ref{def:graph}.

The values in the theorem
are indicative, having been chosen for their
readability. More precise values of $c(d)$
for each $d$ can be obtained, 
but they are of limited interest.
More interesting is to 
mention that $c(d)\rightarrow \frac 4
3$ when $d$ becomes large.
These results 
on $c(d)$ provide  significantly better  lower bounds for $\gamma_{d+1}$
that was previously known $\gamma_{d+1}\ge c(d)$ improving upon $\gamma_{d+1} \ge 0.48$
in the case~2, 
and improving upon $\gamma_{d+1} \ge \frac 2 3$ in the case~3.
The fact that $c(d)\le \frac 4 3$ shows that no further
improvement can be expected from the trick introduced in the
present paper.

Furthermore, these explicit families of graphs
do even better than what
the probabilistic method~\cite{ErSa63}
is able to achieve, namely  a $\gamma_d\ge 1$.
When dealing with Cayley graphs on
$PGL_2(\F_q)$, it was
proved in Theorem 9 of~\cite{GaHoShShVi09}  that random Cayley graphs\footnote{the model of random Cayley graphs is described
p.~2 of~\cite{GaHoShShVi09}} have a girth $\ge 
( \frac 13  - o(1)) \log_d |PGL_2(\F_q)|$ for $q$ sufficiently large.
The exact value is not known, but the
new graphs of the present paper have most likely much larger girth 
than the one for the corresponding
random Cayley graph.

\paragraph{The main inequality}
This paragraph presents the main intermediate
result~\eqref{eq:girth}, and the next paragraph will show
how to deduce from it  the
bounds of Theorem~\ref{th:contrib}.
A few more notations are necessary:
\begin{Def}\label{def:kappa}
Given an integer $d$, $p$ denotes any prime number
$p\ge d$, with the additional condition $p\equiv 3 \bmod 8$
when $d$ is even.
Let $\kappa:= \log_p d\ge 1$,
so that $p = d^\kappa$. Define $$Q_d(p):=\max\{ p^8 , 120^\kappa p\}.$$
\end{Def}
Given another prime $q > Q_d(p)$,
there is a symmetric\footnote{that is if $x\in \Dd_{p,q}$,
then $x^{-1} \in \Dd_{p,q}$ as well}
 subset $\Dd_{p,q}$ of $PGL_2(\F_q)$
of cardinality $d+1$, such that if we define:
\begin{equation*}
G_{d,p,q}:=\begin{cases}
\Cayl(PGL_2(\F_q),\Dd_{p,q}) & \text{if $\Jacobi{p}{q}=-1$}\\
\Cayl(PSL_2(\F_q),\Dd_{p,q}) & \text{if $\Jacobi{p}{q}=1$}
\end{cases}
\end{equation*}
(See Definition~\ref{def:graph}
for more details on $G_{d,p,q}$). Then:
\begin{itemize}
\item $G_{d,p,q}$ is a $d+1$-regular graph of size
  $|PGL_2(\F_q)|=q^3-q$ or $|PSL_2(\F_q)|=\frac 1 2 (q^3-q)$
according to the sign of the
Legendre symbol $\Jacobi{p}{q}$.
\item $G_{d,p,q}$ is connected,  bipartite if
$\Jacobi{p}{q}=-1$,
  and not bipartite if $\Jacobi{p}{q}=1$.
\item the girth of $G_{d,p,q}$ satisfies the {\bf main inequality}:
\begin{equation}\label{eq:girth}
\girth (G_{d,p,q}) \ge
\begin{cases}
\frac{2}{3\kappa} \log_d |G_{d,p,q}| & \text{if $\Jacobi{p}{q}=1$}\\
\frac 4 {3\kappa}\log_d |G_{d,p,q}| - \log_p 4 & \text{if $\Jacobi{p}{q}=-1$}
\end{cases}
\end{equation}
\end{itemize}
Let us point out here
that $\girth(G_{d,p,q})\le \frac 4 3 \log_d |G_{d,p,q}| + 1$
or $\girth(G_{d,p,q})\le \frac 2 3 \log_d|G_{d,p,q}| + 1$,
for any $d$. Indeed, these lower bounds already occur for the
Ramanujan graphs~\cite[Last proposition]{Ma88},
from which the graphs
$G_{d,p,q}$ are derived. This is why $c(d)\le \frac 4 3$ in Theorem~\ref{th:contrib}.

Fixing $p$ and $d$, we can consider the following two kinds of infinite families
of graphs, indexed by $q$:
\begin{equation}\label{eq:X}
\X_{d,p}:=\{G_{d,p,q}\}_{q\ \text{prime},\ q>Q_d(p),\ \Jacobi{p}{q}=-1},
\end{equation}
and
\begin{equation}\label{eq:Y}
\Y_{d,p}:=\{G_{d,p,q}\}_{q\ \text{prime},\ q>Q_d(p),\ \Jacobi{p}{q}=1}.
\end{equation}
From Main Equality~\eqref{eq:girth} above, we infer: $\gamma(\X_{d,p})\ge \frac 4
{3\kappa}$
and $\gamma(\Y_{d,p})\ge \frac 2 {3\kappa}$, where $\kappa = \log_d p$.

\paragraph{Main Inequality implies Theorem~\ref{th:contrib}}
It is quite easy to recover the bounds on $c(d)$ of
Theorem~\ref{th:contrib}
from Main Inequality~\eqref{eq:girth}.
The lower bound on the girth in~\eqref{eq:girth}
is indeed the largest when $\kappa$ is the smallest.
To minimize $\kappa$, let us first introduce some notations:
\begin{Def}\label{def:graph}
Given an integer $u>5$, let
\begin{equation*}
p(u):=\min \{p\ge u \ : \ p\ \text{prime}\},
\quad \text{and}
\quad p_3(u):= \min \{p\ge u\ :\ p\ \text{prime $\equiv 3
\mod 8$}\}.
\end{equation*}
Then, for each $d\ge 10$, we consider two families
of graphs $\X_d$ and $\Y_d$ as:
\begin{multline*}
\text{if $d$ is even:}\quad \X_d:=\X_{d,p(d)},
\qquad \Y_d:=\Y_{d,p(d)},\\
\qquad \text{and if $d$ is odd:}
\quad \X_d:=\X_{d,p_3(d)},
\qquad \Y_d:=\Y_{d,p_3(d)}
\end{multline*} 
The real number $\kappa$
of Definition~\ref{def:kappa}
verifies then $\kappa = \log_d p(d)$
if $d$ is odd
and, $\kappa = \log_d p_3(d)$ if $d$ is even.
\end{Def}
Then, minimizing $\kappa$ brings in the question: Given $u$ odd,  how big is the smallest prime $p(u)$
larger than $u$ ?
Similarly , if $u$ is even, how big can $p_3(u)$ be?
\smallskip

Considering the worst case where $u$ is equal to a prime plus one,
this is related to the problem of {\em gap between primes}~\cite[pp.10-12]{Gr95}.
Bertrand's postulate affirms that $p(u)<2u$, Cram\'er's conjecture suggests that
$p(u)<\log(u)^2$ for some reasonably large $u$, in between various upper bounds on the gap between two primes
have appeared, most being valid only for ``large enough'' values of $u$.
For us, small values of $u$
must be taken into account and therefore we use the unconditional estimate
$p(u)<u(1+\frac{1}{2 (\log  u)^2})$ valid for $u\ge 3275$ (see~\cite[Sec.~4]{Du99}).
Sharper estimates would yield (tiny) better estimates on the girth
only for large degrees of regularity $d$, when $1.333 < c(d)<4/3$.

It implies that: $\kappa\le \log_{u} (u(1+\frac{1}{2 (\log u)^2}))$
for $u\ge 3275$, and proves that $c(d)=\frac{4}{3 \kappa}\ge 1,33$
for $d\ge 3275$. 
For smaller values of $d$, I used a computer
and found the following.
The smallest integer $d_1$ for which
$\big[d\ge d_1
\Rightarrow$ $\frac{4}{3\kappa}\ge 1,33$ 
with $\kappa=\log_{d_1} p(d_1)\big]$ is 1335, and then $p(1335)=1361$.
The smallest integer $d_2$ 
for which
$\big[d\ge d_2\Rightarrow
\frac{4}{3\kappa}\ge 1,3$
with $\kappa=\log_{d_1} p(d_1)\big]$ is 35, and then $p(35)=37$.
Between 15 and 31, it is easy to check that $\frac 4 {3\kappa} \ge
1,27$.
There is no integer smaller than 15 and  greater than 10 which is not
a prime power.
This achieves the  proof of the bound on $c(d)$ in
Theorem~\ref{th:contrib},
when $d$ is odd.

As for $p_3(u)$, I used results of~\cite{RaRu96}.
This requires to introduce the classical arithmetic function
$$
\theta(x;k,\ell):=\sum_{p\equiv \ell\bmod k\atop p\le x} \ln(x),
\qquad
\text{where $p$ denotes a prime number}.
$$
Indeed, there is  a prime number
equal to 3 modulo 8 in the interval $[a;b]$
if $\theta(b;8,3) - \theta (a;8,3) >0$.
The estimate of~\cite[Theorem~1]{RaRu96} shows:
$$
\max_{1\le y\le x} |\theta(y;8,3)-\frac{y}4|\le 0,002811 \frac x
  4,\qquad \text{for $x\ge 10^{10}$}.
$$
Setting $\epsilon=0,002811$, for $x\ge 10^{10}$ and any $y$, it comes: 
$$
\frac y  4 -\epsilon \frac x 4 \le \theta(y;8,3) \le \epsilon \frac x 4
+\frac y 4.
$$
It follows that for all $b>a\ge 10^{10}$, 
$$
\theta(b;8,3) - \theta(a;8,3) \ge \frac b 4 (1 - 2\epsilon) - \frac a 4.
$$
This insures that for $a\ge 10^{10}$ there is a prime equal to 3
modulo 8 in each interval $[a;\frac a{1-2 \epsilon}]$.
For $d\ge 10^{10}$, this clearly proves that
$\frac 4{3\kappa}\ge 1,33$, since then $\kappa=\log_d p_3(d) \le \log_d
\frac{d}{1-2 \epsilon}$. 
For values $d\le 10^{10}$, a laptop computer may not be powerful enough to
check what the maximal values of $\log_d p_3(d)$ are.
Again, from~\cite[Theorem~2]{RaRu96}, in this case:
$$
\max_{1\le y \le x} | \theta(y;8,3) - \frac y 4| \le  1,82 \sqrt{x},
\qquad
\text{for $1 \le x \le 10^{10}$}.
$$
It follows that $\theta(b;8,3) - \theta(a;8,3) \ge \frac{b-a}4 - 2\cdot
1,82 \sqrt{b}$ for $b > a$. This shows that in the interval
$[a;a(1+\frac{8\cdot 1,82}{\sqrt{a} - 8\cdot 1,82})]$ there is
a prime equal to 3 modulo 8. Hence, $\kappa = \log_d p_3(d)\le 1 +
\log_d (1+\frac{8\cdot 1,82}{\sqrt{d} - 8\cdot 1,82})$,
showing that $\frac 4 {3 \kappa}\ge 1,33$ if $d\ge 228050$.

The other values of $c(d)$ of Theorem~\ref{th:contrib}
in the case $d$ even,  for $ d \le 228050$
are easily obtained with the help of a computer.
This concludes the proof of Theorem~\ref{th:contrib}
assuming the main inequality~\eqref{eq:girth}.

\section{Proof of the main inequality}
It remains to show that Main Equality~\eqref{eq:girth} holds.
All the necessary material is contained in the 
monograph~\cite{DaSaVa03}. To make this section
a minimum self-contained, many results appearing therein  are
recalled.

\subsection{Unique factorization of quaternions and regular trees\label{sec:unique}}
The construction of Ramanujan graphs by
Lubotzky-Philips-Sarnak is achieved by taking finite quotients
of a ``mother graph'', which is  a regular tree.
They used  simply
the factorization of quaternions to build these regular trees.

We briefly recall  this here, referring to Ch.~2.6 of the
aforementioned monograph~\cite{DaSaVa03} for the details.

\paragraph{Quaternions}
For $R$ a commutative ring,
let $\H(R)$ denotes the Hamilton quaternion algebra over $R$:
$$
\H(R):= R + R\ii + R \jj + R \kk,\qquad \ii^2=\jj^2=\kk^2=-1,\quad \kk=
\ii \jj=-\jj\ii.
$$
The {\em conjugate} of an element $\alpha = a_0 + a_1\ii + a_2\jj +a_3\kk$
is $\conj{\alpha}:=2 a_0 - \alpha = a_0 -a_1 \ii - a_2\jj -a_3 \kk$,
and the {\em norm} of $\alpha$ is $N(\alpha)=\alpha
\conj{\alpha}=a_0^2+a_1^2+a_2^2+a_3^2$.
\smallskip
The multiplication of quaternions makes the norm  multiplicative: $N(\alpha \beta)=N(\alpha)N(\beta)$.
\smallskip
Given a quaternion $\alpha=a_0+a_1\ii + a_2\jj + a_3 \kk$
the non-negative integer $\gcd(a_0,a_1,a_2,a_3)$
is called the {\em content} of $\alpha$ and is denoted $c(\alpha)$.
If $c(\alpha)=1$, then $\alpha$ is {\em primitive}.
\smallskip

Let us set  $R=\Z$. We introduce a property of unique
factorization
for integral quaternions $\H(\Z)$, yet in a special easy case that is
sufficient
for the purpose of this article. This restriction
is to consider only  quaternions
whose norm  is a power of an odd prime $p$ (instead of considering any
quaternion in $\H(\Z)$).
\smallskip

Given an odd prime $p$, and a primitive  quaternion $\alpha\in \H(\Z)$
of norm $p^k$, then there exist  {\em prime quaternions}
$\pi_1 , \ldots , \pi_k$  ({\em prime} means that if $\pi=\gamma \delta$,
then either $\gamma$ or $\delta$ is a unit in $\H(\Z)$) such that
$\alpha = \pi_1 \cdots \pi_k$.
In  a word, this follows from the possibility to perform a Euclidean
division in $\H(\Z)$ of two such quaternions whose norm is a power of $p$;
a non-commutative  Euclidean
algorithm (one ``one the right'', one ``on the left'' ) is deduced,
in order to compute left and right gcds.
This permits to show that prime quaternions are precisely
those whose norm is a prime number.
Then the existence of a factorization follows easily by induction
on the exponent $k$  of the norm $p^k=N(\alpha)$.
\smallskip

The default of uniqueness is completely related to the units of $\H(\Z)$
(which are $\pm 1,\pm \ii,\pm \jj , \pm \kk$).
What this means is that two distinct factorizations
$\pi_1 \cdots \pi_k$ and $\mu_1 \cdots \mu _k$
of $\alpha$ verify: $\pi_i=\epsilon_i \mu_i$,
for some $\epsilon_i\in\H(\Z)^\star$ and for $1\le i \le k$.
The group of 8 units $\H(\Z)^\star$ acts
on the set of quaternions of norm $p$.
By isolating one quaternion per orbit,
uniqueness can be recovered.
Since the number of quaternions of norm $p$
is $8(p+1)$ by a famous theorem of Jacobi
(indeed, such quaternions $x_0+x_1\ii+x_2\jj+x_3\kk$
give a solution in $\Z^4$
of $f(x)=p$, where $x=(x_0,x_1,x_2,x_3)$
and $f(x)=x_0^2+x_1^2+x_2^2+x_3^2$).
As perfectly explained in p.~67-68 of~\cite{DaSaVa03},
a quite natural way to isolate one quaternion per orbit is to introduce:
\begin{equation}\label{eq:Pp(p)odd0}
\Pp(p)=\{ \pi \in \H(\Z)\ \text{primitive} : \ N(\pi)=p,\ \pi_0>0,\ \pi-1 \in
2\H(\Z)\} \qquad \text{if $p\equiv 1\mod 4$},
\end{equation}
\begin{multline}
\label{eq:Pp(p)odd}\Pp(p)=\{ \pi \in \H(\Z)\ \text{primitive} : \ N(\pi)=p,\ \ \pi_0 >0
\ \ \text{if}
\ \ \pi_0\not=0,\ \  \text{ and } \ \   \pi_1>0 \ \ \text{ else,}\\
\pi - \ii - \jj -\kk \in  2\H(\Z)\} \qquad \text{if $p\equiv 3\mod 4$}
\end{multline}
The fact that $\pi_0\not=0$,
or $\pi_1\not=0$ if $\pi_0=0$ is made clear by the explanations coming
hereafter.
\begin{Rem}\label{rem:Pp(p)}
Some general remarks about this set:
\begin{enumerate}[(a)]
\item \label{rem:1} if $\alpha \in \Pp(p)$, then $\epsilon \alpha$
and $\alpha\epsilon$ are not in $\Pp(p)$,
 for any unit $\epsilon \in \H(\Z)^\star$ different from 1.
\item \label{rem:2} Similarly, given $\beta\in\H(\Z)$,
$N(\beta)=p$, there are exactly two units $\varepsilon,\varepsilon' \in \H(\Z)^\star$
that yield $\varepsilon \beta\in \Pp(p)$
and $\beta \varepsilon'\in\Pp(p)$.
\item this implies that $|\Pp(p)|=p+1$ (according to Jacobi's theorem on the sum
of four squares).
\item \label{rem:4}given $\pi\in \Pp(p)$, if $\pi_0\not=0$ then $\conj{\pi}\in \Pp(p)$
(easy to check). If $\pi$ is such that $\pi_0=0$, as it may happen
when $p\equiv 3 \bmod 4$ (actually when $p\equiv 3 \bmod 8$ after Proposition~\ref{prop:38},
then $\conj{\pi} = -\pi \not\in \Pp(p)$,
in conformity  with the two points~\eqref{rem:1} and~\eqref{rem:2} above.
\end{enumerate}
\end{Rem}
Remark that the first point~\eqref{rem:1} allows
a form of uniqueness of the factorization of
quaternions~\cite[2.6.13~Theorem]{DaSaVa03}.
\begin{Theo}\label{th:fact}
Given $\alpha$ of norm $p^k$, and of content $c(\alpha) = p^\ell$, then
there exist unique $\pi_1,\ldots,\pi_{k-2\ell}\in \Pp(p)$
and a  unique unit $\epsilon \in \H(\Z)^\star$ such that:
\smallskip

\centerline{$\alpha =  c(\alpha) \, \epsilon\, \pi_1 \cdots \pi_{k - 2 \ell},
\qquad \text{with}\ \ \conj{\pi_i}\not= \pi_{i-1} \ \text{ if $\ \conj{\pi_i}
\in \Pp(p),\ $ and with  $\ \pi_i \not= \pi_{i-1}\ $ else.}$\foorp}
\end{Theo}
Let us stress that under these conditions,
the quaternion $\pi_1 \cdots \pi_{k - 2 \ell}$ is primitive
(motivating the definition of
irreducible product in Definition~\eqref{def:irredProd}).
\medskip

We focus on the case $\pi\in \Pp(p)$ and $\conj{\pi}\not\in \Pp(p)$,
which may happen when $p\equiv 3 \bmod 4$ as mentioned
in~\eqref{rem:4} above.
\begin{Prop}\label{prop:38}
There is an element
$\pi = \pi_0 + \pi_1 \ii + \pi_2\jj + \pi_3\kk\in \Pp(p)$
for which $\pi_0=0$ (equivalently
$\pi=\conj{\pi}$, or $\conj{\pi}\not\in \Pp(p)$)
 if and only if $p \equiv 3 \bmod 8$.
\end{Prop}
\proof
By definition of $\Pp(p)$
this can only happen if  $p\equiv 3 \bmod 4$,
since otherwise $\pi_0\equiv 1\bmod 2$.
For such a $\pi$, $N(\pi) = \pi_1^2+\pi_2^2+\pi_3^2$
and consequently  $p$ is a sum of 3 squares.
Reciprocally, a sum of 3 squares $x_1^2+x_2^2+x_3^2$ equal to  $p$
gives a quaternion $x = x_1\ii + x_2\jj + x_3\kk \in \H(\Z)$
of norm $p$, which is also necessarily primitive (because
$p$ is prime). Since $p\equiv 3\bmod 4$, $p$ is not the sum of 2 squares. Hence, necessarily,
 $x_1\equiv x_2\equiv x_3\equiv 1\bmod 4$, implying $x\in \Pp(p)$.

We have proved that such a $\pi$ exists in $\Pp(p)$ if and only if
$p\equiv 3 \bmod 4$ and $p$ is the sum of 3 squares.
This is true if and only if $p\equiv 3 \bmod 8$, as the Gauss
theorem~\cite[Ch.~IV, {\em Appendice}]{Se77arithm} on  sum of 3 squares shows:
\begin{Theo}[Gauss]
An integer $n$ is the sum of 3 squares if and only if
$n$ is not equal to $4^k(8 \ell +7)$ whatsoever are $k,\ell \in \N$.\foorp
\end{Theo}
Hence, $p\equiv 3 \bmod 4$ is the sum of 3 squares
if and only if $p\not = 4^k(8\ell +7)$.
Suppose $p=4^k(8\ell +7)$, then $p\equiv 3 \bmod 4$
gives $k=0$, and $p=8\ell + 7$, implying $p\equiv 3 \bmod 8$.
This proves that $p\equiv 3\bmod 4$
is sum of 3 squares if and only if $p\equiv 3 \bmod 8$,
achieving the proof of Proposition~\ref{prop:38} \foorp
\medskip

In the case $p\equiv 3 \bmod 8$, we denote a quaternion in $\Pp(p)$
of the form shown in Proposition~\ref{prop:38} by
the letter $\nu$,
and the others by the letter $\mu$. One has:
\begin{equation}\label{eq:Pp(p)38}
\text{if $p\equiv 3 \bmod 8$},\quad
\Pp(p)=\{\mu_1,\ldots , \mu_s,\nu_1,\ldots, \nu_t\},
\quad \text{with $s + t = p + 1$ and $t>0$}.
\end{equation}
Note that $s$ is even because each $\mu_i$ comes along with its conjugate,
that is there is an $i'\not=i$ such that $\conj{\mu_i}=\mu_{i'}$. That is,
$t=p+1-s$ is also even.

\paragraph{Trees built on quaternions}
The unique factorization theorem~\ref{th:fact}
permits to build infinite regular trees of arbitrary degree $d$.
As in Definition~\ref{def:kappa}, let $p\ge d$ be a prime number, ordinary if $d$ is odd,
and equal to 3 modulo 8 if $d$ is even.

\begin{Lem}\label{lem:Dp}
We can choose a subset $\Dd(d)\subset \Pp(p)$ of cardinality $d+1$
such that, given  $\pi \in \Dd(d)$,
one has: $\conj{\pi}\in \Dd(d)$ if  and only if $\conj{\pi}\in \Pp(p)$.

In particular, if $d$ is even then $\Dd(d)$ contains at least
one $\pi$ such that $\conj{\pi} \not\in \Dd(d)$ (this latter
case happens only if $p\equiv 3 \bmod 8$
according to Proposition~\ref{prop:38}).
\end{Lem}
\proof
If $d$ is odd, then Property~\eqref{eq:Pp(p)odd0} in the definition~\eqref{eq:Pp(p)odd} of $\Pp(p)$
when $p\not\equiv 3 \bmod 8$ makes it clear: it suffices
to choose $\frac{d+1}2$ elements pairwise not conjugate,
as well as  their $\frac{d+1}2$ conjugates (that are
also in $\Pp(p)$ in this case).
For the case $p \equiv 3 \bmod 8$, let us use the two
even integers $s$ and $t$ defined in~\eqref{eq:Pp(p)38}.
We first choose $k_1:=\max\{ \frac{d+1}2 , \frac s 2\}$ couple of conjugates
in $\Pp(p)$, and, if necessary, $d+1- 2 k_1$ elements
$\pi$ such that $\conj{\pi}\not\in \Pp(p)$.

If $d$ is even, then $p\equiv 3\bmod 8$ by
Definition~\ref{def:kappa}.
A way to choose the set
$\Dd(d)$ is as follows. First choose $k_1 := \max\{ \frac d 2  ,
\frac s 2\}$
couples of conjugates, completed with 
$d + 1 - 2 k_1$ elements $\pi$ such that $\conj{\pi}\not\in\Pp(p)$.
\foorp

Notice that in general, there are several other possible ways of
choosing $\Dd(d)$ inside $\Pp(p)$.
\begin{Def}\label{def:irredProd}
An {\em irreducible product} of length $\ell$
over $\Dd(d)$ is the product of $\ell$ elements
$\alpha_1,\ldots,\alpha_\ell$
in $\Dd(d)$ where two consecutive elements:

-  are not conjugate, $\alpha_i \not= \conj{\alpha_{i+1}}$, if
$\conj{\alpha_i} \in \Pp(p)$

-  are not equal, $\alpha_i\not= \alpha_{i+1}$, if $\conj{\alpha_i}\not\in\Pp(p)$. 

\noindent The set of all irreducible products over $\Dd(d)$ is denoted
$\Lambda_\Dd$.
\end{Def}
The motivation of this terminology comes from
the following fact, resulting of the unique factorization (Theorem~\ref{th:fact}):
the product of a sequence of elements in $\Dd(d)$
that does not verify the conditions mentioned in the definition
can be reduced, yielding a non primitive quaternion. 

Furthermore,
Theorem~\ref{th:fact} also tells that  two different irreducible products yields
two different quaternions. This allows to define a $d+1$-regular
tree $T_{d}$ in the following way:
\begin{itemize}
\item the vertex set $V(T_d)$ is identified with the irreducible products of
$\Lambda_\Dd$ over $\Dd(d) \subset \Pp(p)$
\item the root is identified with the void product;
given another vertex identified with the irreducible
product $\alpha_1 \cdots \alpha_s$, we define $d$ adjacent
vertices whose  irreducible products are:
$$
\alpha_1\cdots \alpha_s \alpha_{s+1},\qquad
\alpha_{s+1}\in \Dd(d) \quad \text{where}
\begin{cases}
\alpha_{s+1} \not= \conj{\alpha_s} & \text{if $\conj{\alpha_s} \in \Pp(p)$}\\
\alpha_{s+1} \not= \alpha_s & \text{if $\conj{\alpha_s} \not\in \Pp(p)$}
\end{cases}
$$
\item and the last adjacent vertex is the irreducible product
$\alpha_1\cdots \alpha_{s-1}$
\end{itemize}

\subsection{Algebraic construction of the tree
and definition of the graphs $G_{d,p,q}$\label{sec:Gdpq}}
It is necessary to give an interpretation of the tree
$T_d$ constructed above more algebraically.
Indeed, the graphs $G_{d,p,q}$ are naturally defined algebraically.
\paragraph{Algebraic construction of the trees $T_{d}$}
It consists in seeing the trees $T_d$ as Cayley graphs on free groups if $d$ is even,
or on groups with involutions and unique factorization
property in terms of the generating set if $d$ is odd.
These groups are: 
\begin{Prop}
The set  $\Lambda_\Dd$ of all irreducible
products over $\Dd(d)$ can be endowed
with the structure of group. Id $d$ is odd, then $\Lambda_\Dd$ is free over $\Dd(d)$.

If $d$ is even, $\Dd(d)$ contains at least one involution --- hence $\Lambda_\Dd$ is not free on $\Dd(d)$ ---
but each element of $\Lambda_\Dd$ can be uniquely
written as a product of elements in $\Dd(d)$.
\end{Prop}
\proof
Given two irreducible products $\alpha := \alpha_1 \cdots \alpha_n$,
and $\beta:=\beta_1 \cdots \beta_m$ in $\Lambda_\Dd$,
we associate an irreducible product $\alpha \times \beta $ as follows.

-   there is no integer $i\ge 0$ such that
$\alpha_{n-i}
\not=\conj{\beta_{i+1}}$ if $\conj{\beta_{i+1}} \in \Pp(p)$,
or $\alpha_{n-i} \not = \beta_{i+1}$ if $\conj{\beta_{i+1}}\not\in \Pp(p)$.
Then we define $\alpha \times \beta =1$.

- else, let $\ell \ge 0$ be the largest such integer $i$.
Then the content of $\alpha \beta$ is then $c(\alpha\beta)=p^\ell$,
and $\frac{\alpha \beta} {p^\ell}$ is primitive.
Its unique
factorization is  given by:
$\frac{\alpha \beta} {p^\ell}=
\pm \alpha_1\cdots \alpha_{n-\ell} \beta_{\ell+1} \cdots \beta_m$.
This allows to define,
$$
\alpha \times \beta := \alpha_1\cdots \alpha_{n-\ell} \beta_{\ell+1} \cdots
\beta_m.
$$
Note that this is an irreducible product in $\Lambda_\Dd$.
\smallskip

It is easy to check that $\times$ defines
an associative operations
on $\Lambda_\Dd$ with unit element 1
(the void irreducible product).
The inverse of an irreducible product
$\alpha:=\alpha_1\cdots \alpha_n$ is $\beta:=\tilde{\alpha_n}\cdots
\tilde{\alpha_1}$ where $\tilde{\alpha_i} = \conj{\alpha_i}\in \Dd(d)$
if $\conj{\alpha_i}\in \Pp(p)$,
and $\tilde{\alpha_i}=\alpha_i$ if $\conj{\alpha_i}=-\alpha_i\not\in \Pp(p)$.
The content of $\alpha\beta$ is then $p^n$, hence $\alpha \times \beta = 1$.
\smallskip

It remains to show that each element of the group $(\Lambda_\Dd, \times)$ 
can be uniquely written as a product of elements in $\Dd(d)$.
This follows by the  definition~\ref{def:irredProd}
of irreducible products
on $\Dd(d)$, that  yields different quaternions by the
unique factorization theorem~\ref{th:fact}.
\foorp

\begin{Rem}\label{rem:free-group}
Using the notations in~\eqref{eq:Pp(p)38}, 
$\Dd(d)$ consists of elements
$\mu_1,\ldots, \mu_i,\nu_1,\ldots,\nu_v$
with $u\le s$ and $v \le t$,
such that $\conj{\nu_i}\notin \Dd(d)$ and $\conj{\mu_i}=\mu_{i'}\in \Dd(d)$.
Let $(K,\times)$ be the subgroup of
$(\Lambda_{\Dd},\times)$ 
generated by $\mu_1,\ldots,\mu_u$. 
This is a free group for $\times$,  and we have:
$$
(\Lambda_\Dd , \times ) \simeq (K,\times) * \l \nu_1\r \cdots *\l \nu_v\r,
$$
where $\l \nu_i\r$ is the subgroup of order 2 of
$(\Lambda_\Dd ,\times)$
 generated by $\nu_i$  and $*$  is the free product on subgroups
of $(\Lambda_\Dd , \times )$.\end{Rem}

The combinatorial definition of the tree $T_d$ given
at the end of Section~\ref{sec:Gdpq} is the Cayley
graph of the group $\Lambda_\Dd$ with generating set
$\Dd(d)$.
\begin{equation*}
T_d \simeq \Cayl(\Lambda_\Dd , \Dd(d)).
\end{equation*}

\paragraph{Graphs $G_{d,p,q}$ as finite quotients of the tree $T_d$} 
As above, we let $d$ be an integer greater than 10,
and $p$ a prime greater than $d$, equal to 3 modulo 8 id $d$ is even
(and without condition if $d$ is odd).
Now we let $q>Q_d(p)$ where $Q_d(p)$ is the constant introduced
in Definition~\ref{def:kappa}.

The next step consists in taking finite quotients of
the tree $T_d$.
Let 
\begin{equation}\label{eq:tau_q}
\tau_q:\H(\Z) \rightarrow \H(\F_q)
\end{equation}
the reduction map modulo $q$.
When restricted to $\Lambda_{\Dd}$,
we observe the following:
\begin{itemize}
\item $\tau_q(\Lambda_{\Dd} ) \subset \H(\F_q)^\star$
\item  $\tau_q(\alpha \beta)$ and $\tau_q(\alpha * \beta)$ differ multiplicatively
by $\tau_q(p^{\ell})$, where $p^\ell$ is the content of $\alpha \beta$,
which is in the center $\Zz$ of the group $\H(\F_q)^\star$.
\end{itemize}
Hence, by taking the quotient group $\H(\F_q)^\star/\Zz$
the following map:
$$
\mu_q : \Lambda_{\Dd} \rightarrow \H(\F_q)^\star/\Zz,
$$
is a group homomorphism.
Next, we  identify the image of this group homomorphism.
Recall that since $p\not=2$,
the quaternion algebra over $\F_q$ as defined in Section
~\ref{sec:unique}
is isomorphic to the algebra of 2-by-2 matrices over $\F_q$.
Indeed, in $\F_q$ there are two elements
$x$ and $y$ such that $x^2+y^2+1=0$ (see Prop.~2.5.2 and~2.5.3
in~\cite{DaSaVa03}). The following map is an isomorphism of $\F_q$-algebra:
\begin{eqnarray*}
\phi \ :\ \ \H(\F_q) & \rightarrow & M_2(\F_q),\\
  \alpha_0 + \alpha_1 \ii + \alpha_2 \jj + \alpha_3\kk & \mapsto &
\left(\begin{array}{cc} \alpha_0 +\alpha_1 x + \alpha_3 y & -\alpha_1
  y +\alpha_2 +\alpha_3 x\\
-\alpha_1 y -\alpha_2 +\alpha_3 x & \alpha_0 -\alpha_1 x -\alpha_3 y
\end{array}\right).
\end{eqnarray*}
Moreover $N(\alpha)=\det \phi(\alpha)$. 
We deduce the following group isomorphism $\psi$
from $\phi$:
\[\label{eq:psi}
\psi:\H(\F_q)^\star/\Zz \rightarrow PGL_2(\F_q),
\]
and we let:
$$
\conj{\mu_q}:=\psi \mu_q ,\quad \text{and}
\quad  \ker \conj{\mu_q} := \Lambda_\Dd(q),
\quad \text{so that} \quad \Lambda_\Dd /\Lambda_\Dd(q) \hookrightarrow
PGL_2(\F_q).
$$
\begin{Lem}\label{lem:ps12}
If $q$ is such that $p$ is a quadratic residue modulo $q$,
then $\conj{\mu_q}(\Dd(d)) \subset PSL_2(\F_q)$. Else,
$\conj{\mu_q}(\Dd(d)) \subset
PGL_2(\F_q) -  PSL_2(\F_q)$.
\end{Lem}
\proof 
The group homomorphism
$\epsilon:\H(\F_q)^\star \rightarrow \{-1,1\}$, $x\mapsto
\Jacobi{N(\alpha)}{q}$
takes the same value
on each class modulo the center $\Zz$.
The factor map   $\conj{\epsilon} : \H(\F_q)^\star / \Zz \rightarrow
\{-1,1\}$, $x\Zz \mapsto \epsilon(x)$,
is well-defined.
The set of quaternions in $\H(\F_q)^\star$
of norm 1, denoted $\H_1$,  is sent to 1 by $\epsilon$,
and hence $\ker \conj{\epsilon} \supset \H_1/(\Zz\cap \H_1) $.
Now, given $\pi\in\Dd(d)$,  
$\Jacobi{p}{q}$ and $\conj{\epsilon} (\mu_q(\pi)) $
are equal. This shows that if $\Jacobi{p}{q}=1$, then
$\mu_q(\Dd(d))\subset \ker \conj{\epsilon}$, and if
$\Jacobi{p}{q}=-1$, then $\mu_q(\Dd(d)) \subset \H(\F_q)^\star/\Zz  - 
\ker \conj{\epsilon}$. Using the isomorphism $\psi$,
we obtain $\conj{\mu_q}(\Dd(d)) \subset PSL_2(\F_q)$
if $\Jacobi{p}{q}=1$, and $\conj{\mu_q}(\Dd(d)) \subset PGL_2(\F_q)  - 
PSL_2(\F_q)$ else.\foorp
\medskip

By  the above discussion, comes:
\begin{equation}\label{eq:embed}
\Lambda_{\Dd}/\Lambda_{\Dd}(q) \hookrightarrow \begin{cases}
 PSL_2(\F_q) & \text{if $\Jacobi{p}{q}=1$}\\ 
 PGL_2(\F_q) & \text{if $\Jacobi{p}{q}=-1$}\end{cases}
\end{equation}
\smallskip

\begin{Lem}\label{lem:Dpq}
Let $\Dd_{p,q}:=\conj{\mu_q} (\Dd(d))$.
One has $|\Dd_{p,q}|=|\Dd(d)|= d+1$
\end{Lem}
\proof
The map $\psi$~\eqref{eq:psi} being an isomorphism it suffices to show that
$|\Dd(d)|=|\mu_q(\Dd(d))|$. Since $\Dd(d)\subset\Pp(p)$,
this will certainly follow from $|\Pp(p)|=|\mu_q(\Pp(p))|$.
The later is (easily) proved in~\cite[4.2.1~Lemma]{DaSaVa03},
under the assumption 
that $q>2\sqrt{2}$, verified because $q>Q_d(p)\ge p^8$.
\foorp

Already mentioned in the Introduction, we can now
give a precise definition of the graph  $G_{d,p,q}$:
\begin{Def}\label{def:Gdpq}
Given  the three integers $d$, $p$
and $q$ as defined above, the graph
$G_{d,p,q} $ is :
$$
G_{d,p,q}:=\begin{cases}
\Cayl(PGL_2(\F_q) \, ,\, \Dd_{p,q}) & \text{if $\Jacobi{p}{q}=-1$}\\
\Cayl(PSL_2(\F_q)\, ,\, \Dd_{p,q})   & \text{if $\Jacobi{p}{q}=1$}\end{cases}$$
\end{Def}
By Lemma~\ref{lem:Dpq}, the graphs $G_{d,p,q}$ are $d+1$-regular. Moreover:
\begin{Lem}\label{lem:bi}
The graphs
$G_{d,p,q}$ are bipartite when $\Jacobi{p}{q}=-1$.

Moreover, assuming that $G_{d,p,q}$ is connected when
$\Jacobi{p}{q}=1$, $G_{d,p,q}$ is non-bipartite.
\end{Lem}
\proof
In the first case, a bipartition ${\cal A} \cup {\cal B}$ of the
set of vertices $V(G_{d,p,q})$ is given by ${\cal A}:= PSL_2(\F_q)$,
and ${\cal B} := PGL_2(\F_q) -  PSL_2(\F_q)$.
Indeed, the index of $PSL_2(\F_q)$ in $PGL_2(\F_q)$
is 2, and
Lemma~\ref{lem:ps12} shows that the 
generating set $\Dd_{p,q}$ lies in 
the non-trivial coset $\subset {\cal B}$.

As for the case $\Jacobi{p}{q}=1$,
saying that $G_{d,p,q}$ is connected is equivalent
to saying that $\Dd_{p,q}$ generates $PSL_2(\F_q)$.
Then a bipartition would imply a non-trivial
group homomorphism $PSL_2(\F_q) \rightarrow \{-1,1\}$,
whose kernel would be a proper normal subgroup
of $PSL_2(\F_q)$, excluded since $PSL_2(\F_q)$ is
simple~\cite[3.2.2~Theorem]{DaSaVa03}.
\foorp

To end this subsection,
let us mention that  all these Cayley graphs are connected
(this is Proposition~\ref{prop:connect},
in particular,  $G_{d,p,q}$ is non-bipartite when $\Jacobi{p}{q}=1$
by the lemma just above).
This point is important for estimating the girth, and is not trivial.
In~\cite{LPS88} the authors resort a deep and technical result
of Malyshev on the number of integer solutions of quadratic definite
positive forms;
the construction of Margulis~\cite{Ma88}
differs slightly from the one of~\cite{LPS88},
where
a density argument (strong approximation
theorem) was used.
In our modified construction of graphs, the connectedness is also
crucial, 
but none of these two proofs would work. 
Fortunately, later appeared in~\cite{DaSaVa03}
(see discussion p.~6)
a simple proof of the connectedness,
based on the properties of the subgroups of $PSL_2(\F_q)$,
observed by Frobenius. It will be
instrumental in the present work.

\subsection{Connectedness and final proof}
Following the method of Ch.~4.3 in \cite{DaSaVa03},
this is achieved by showing logarithmic girth.

Let $X$ denote the connected
component of $G_{d,p,q}$ containing the identity.
\begin{Lem}\label{lem:X=Cayl}
Let $\Dd'(d)$ denotes the image of $\Dd(d)\subset
\Lambda_\Dd$
through the group homomorphism:
$\Lambda_\Dd \rightarrow \Lambda_\Dd/\Lambda_\Dd(q)$.
The following  isomorphism of graphs holds:
$\ \  X\simeq \Cayl(\Lambda_\Dd / \Lambda_\Dd(q) \, ,\,  \Dd'(d))$.  
\end{Lem}
\proof
By definition of Cayley graphs $G_{d,p,q}$, we see that 
$X =\Cayl(\l \Dd_{p,q} \r\, ,\, \Dd_{p,q} )$,
where $\l \Dd_{p,q}\r$ denotes the subgroup of $PGL_2(\F_q)$
generated by $\Dd_{p,q}$.
On the other hand,
since $\Dd(d)$ generates
$\Lambda_\Dd$, $\Dd'(d)$ generates $\Lambda_\Dd / \Lambda_\Dd(d)$.
The embedding~\eqref{eq:embed}
shows that $\Lambda_\Dd/\Lambda_\Dd(q)$ is isomorphic
to a subgroup
of $PGL_2(\F_q)$, which  is precisely $\l \Dd_{p,q}\r$.
This induces the graph isomorphism
$$
\Cayl(\Lambda_\Dd/\Lambda_\Dd(d) \, ,\, \Dd'(d)) \simeq
\Cayl(\l \Dd_{p,q} \r\, ,\, \Dd_{p,q} )
$$
concluding the proof. \foorp

By vertex-transitivity of a Cayley  graph on a group,
the closed paths of length $\ell$
(starting and ending) at a vertex $x$ and the ones (starting and ending)
at a vertex $y$ are in one-one correspondence.
In particular a closed path of minimal length in the graph
is found at each vertex, including the vertex 1.
Thanks to  Lemma~\ref{lem:X=Cayl}, a closed path
starting at the identity of $\Lambda_\Dd / \Lambda_\Dd(q)$
corresponds to a product $\alpha=\alpha_1\times\cdots \times \alpha_t \in \Lambda_\Dd$,
with $\alpha_i \in \Dd(d)$, such that $\alpha\in \Lambda_\Dd(q)$.
Thus:
\begin{equation*}
\girth(X):= \inf \{ t \in \N^\star \ : \ \alpha_1\times\cdots\times\alpha_t\in
\Lambda_{\Dd}(q),\ \ \alpha_i \in \Dd(d)\}. 
\end{equation*}

The computations that follow are classical. They already appeared
in~\cite{LPS88}. 
Note that $x=x_0 + x_1\ii + x_2\jj + x_3\kk \in \Lambda_\Dd(q)$
implies that $q|x_i$ for $i=1,2,3$.
If we write $x_i=q y_i$,
then we see that  $N(x)=x_0^2 + q^2 ( y_1 ^2 + y_2 ^2 + y_3 ^2 ) = p^t$.
At least one $y_i \not = 0$ among the values of $i=1,2,3$,
else $x\not \in \Lambda_\Dd$.
Hence, $t\ge 2 \log_p q = \frac 2 3 \log_q q^3$.

In the case where $\Jacobi{p}{q} = -1$, the graphs $G_{d,p,q}$
are bipartite by Lemma~\ref{lem:bi}
and the girth, as is minimum length of all cycles,
is an even. Hereafter, the girth is equal to $2t$.
A basic refinement is possible in this case:
as before, we get $p^{2t}= x_0^2+q^2(y_1^2 + y_2^2 + y_3^2)$,
with at least one $y_i\not=0$ among $y_1,y_2,y_3$.
Hence, $p^{2t}\equiv x_0^2 \bmod q^2$.
This is equivalent to 
$p^t \equiv \pm x_0 \bmod q^2$,
the group $(\Z/q^2\Z)^\star$ being cyclic.
Therefore, $p^t = \pm x_0 + m q^2$ for  a positive integer $m$.
A simple calculation yields $2 p^t - m q^2 > 0$,
from which $t\ge 2\log_p q -\log_p 2$ follows.
The girth in this case satisfies $\girth(X)\ge \frac 4 3  \log_p q^3 - 2\log_p2$.
\medskip

Recall that $X$ is  the connected component of $G_{d,p,q}$
containing $1$. Its cardinality verifies $|X| \le |PGL_2(\F_q)|
= q^3 - q$, and even $|X|\le |PSL_2 (\F_q)|=\frac 1 2 (q^3-q)$
when $\Jacobi{p}{q}=1$.
The definition~\eqref{def:kappa} of $\kappa$ along with the above show that
$ \frac 2 3 \log_p |X| = \frac {2} {3\kappa} \log_d |X|
\le \frac {2}{3\kappa} \log_d q^3\le
\girth(X)$, if $\Jacobi{p}{q}=1$.
And similarly,
$\frac {4}{3\kappa} \log_d |X| - \log_p 4 \le \girth(X)$
if $\Jacobi{p}{q}=-1$.

The graph $X$ has logarithmic girth.
A trick that  first appeared in~\cite[3.3.4~Theorem]{DaSaVa03}
proves that it implies connectedness. We recall
this theorem resulting from the properties
of subgroups of $SL_2(\F_q)$ due to Frobenius;
a group is said to be {\em metabelian} if
it admits a normal subgroup $N$ such that both
$N$ and $H/N$ are abelian.
It is easy to see that $H$ is metabelian
if and only if for any  four elements
$h_1 , h_2 , h_3 , h_4 \in H$ one has
$$
[[h_1,h_2],[h_3,h_4]]=1,\qquad \text{(where $[a,b]=aba^{-1}b^{-1}$).}
$$

\begin{Theo}[~\cite{DaSaVa03}, 3.3.4~Theorem]
Let $q$ be a prime. Let $H$ be a proper subgroup of
$PSL_2(\F_q)$, such that $|H|>60$. Then $H$ is metabelian.\foorp
\end{Theo}
Hence, to prove that $H=PSL_2(\F_q)$,
it suffices to prove that $|H|>60$ and that $H$ is not metabelian.

\begin{Prop}\label{prop:connect}
Since $d\ge 10$ and $q>\max\{d^{8\kappa},(120 d)^\kappa
\}=\max\{p^8,120^\kappa p\}$, one has that
the graph $G_{d,p,q}$ is connected.
\end{Prop}
\proof
It amounts to show that $X=G_{d,p,q}$.
Thanks to Lemma~\ref{lem:X=Cayl}, it suffices
to show that the embedding~\eqref{eq:embed}
is onto, that is:
$$
\Lambda_{\Dd}/\Lambda_{\Dd}(q) \simeq \begin{cases}
 PSL_2(\F_q) & \text{if $\Jacobi{p}{q}=1$}\\ 
 PGL_2(\F_q) & \text{if $\Jacobi{p}{q}=-1$}\end{cases}
$$
This is equivalent to show that
$\conj{\mu_q} ( \Lambda_\Dd) = PSL_2(\F_q)$ or $PGL_2(\F_q)$.
Since $PSL_2(\F_q)$ is an index 2 normal subgroup
of $PGL_2(\F_q)$ and that $\conj{\mu_q}(\Lambda_\Dd) \not\subset
PSL_2(\F_q)$ if $\Jacobi{p}{q}=-1$,
it suffices to show that
$\conj{\mu_q}(\Lambda_\Dd)\cap PSL_2(\F_q) = PSL_2(\F_q)$.

Let $L:=\conj{\mu_q}(\Lambda_{\Dd}))\cap
PSL_2(\F_q)$.
First, we have $ |L| > 60 $.
Indeed,  by Equation~\eqref{eq:Moore1} and
the bound on the girth of $X$ obtained above,
$$
2 \log_p q \le \girth(X) < 2 \log_d |X| +2,
$$
from which follows $\log_p q - 1 < \log_d |X|$,
then $|X| >  d^{\log_p q -1} = p^{{\frac 1 \kappa} (\log_p q  - 1)}$
and finally $|X| > (\frac q p)^{\frac 1 \kappa}$.

Next, $|X|\le 2 |L|$. The equality may occur if $G_{d,p,q}$ is
connected, i.e. $X=G_{d,p,q}$, and if $\Jacobi{p}{q}=-1$.
It follows that $|L| > \frac  1 2 (\frac q p)^{\frac 1 \kappa}$.
Since, $q\ge 120^\kappa p$, this implies $|L| > 60$.

The second step is to show that $L$ is not metabelian,
that is there exist four elements $\ell_1,\ell_2,\ell_3,\ell_4$
in $L$ such that:
\begin{equation}\label{eq:commutator}
[ [\ell_1,\ell_2],[\ell_3,\ell_4]]\not =1.
\end{equation}
Let 4 elements $\alpha_1,\alpha_2,\alpha_3,\alpha_4$
in $\Dd(d)$.
The commutator $[[\alpha_1,\alpha_2],[\alpha_3,\alpha_4]]$ taken
in the group
$(\Lambda_\Dd , \times)$, yields an
irreducible product of length smaller than 16.
And it is equal to 16 if and only if
$[[\alpha_1,\alpha_2],[\alpha_3,\alpha_4]]$ performed this time
in $\H(\Z)$  is primitive (that is no reduction occurred).

Suppose $\alpha_1,\alpha_2,\alpha_3,\alpha_4$
verifies the latter. Let $\ell_i:=\conj{\mu_q}(\alpha_i)\Lambda_\Dd(q)
\in \Lambda_\Dd/\Lambda_\Dd(q)$.
Then by construction of Cayley graphs,
the commutator $[[\ell_1,\ell_2],[\ell_3,\ell_4]]$
yields 
a backtrackless path
of length 16 in $X$.
Beforehand, we have proved that
$
\girth(X)\ge 2\log_p q
$
which is strictly greater than 16 considering that $q > p^8$.
Hence, we have $[[\ell_1,\ell_2],[\ell_3,\ell_4]] \not = 1$
concluding the proof of~\eqref{eq:commutator}, under the
existence of the $\alpha_i$s in $\Dd(d)$.

It is actually always possible to find such $\alpha_i$s
as soon as $|\Dd(d)|>6$,
as perfectly explained
in the proof of~\cite{DaSaVa03} p.~120, paragraphs (a) and (b). 
This is the case  since $d\ge 10$ by assumption.
\foorp

Since $X=G_{d,p,q}$, it follows that $\girth(G_{p,d,q}) \ge 2 \log_p q
=\frac 2 {3\kappa} \log_d q^3 > \frac 2 {3\kappa} \log_d |G_{p,d,q}|$
if $\Jacobi{p}{q}=1$,
and $\girth(G_{d,p,q})\ge 4 \log_p q - \log_p 4> \frac 4 {3\kappa}
\log_q |G_{d,p,q}| - \log_q 4$ if $\Jacobi{p}{q}=-1$,
achieving the proof of Main Equality~\eqref{eq:girth}.
\smallskip

As for the non-bipartite  $d+1$-regular graphs
$H_n$ mentioned
in Theorem~\ref{th:contrib}, they correspond to
the families $\Y_d$ of Definition~\ref{def:graph}.
It has not be proved yet that they are not bipartite.
Going back to the second point above Main Equality~\eqref{eq:girth},
we must show that $G_{d,p,q}$ is non-bipartite when
$\Jacobi{p}{q}=1$. It was
not possible to prove it at the time of the proof
of Lemma~\ref{lem:bi},
because of the lack of knowledge of the connectedness.
Granted by Proposition~\ref{prop:connect},
this concludes the proof of Theorem~\ref{th:contrib}.


\section*{Concluding remarks}
\paragraph{On the previous work}
By a simple modification made  on the classical construction of Ramanujan graphs
of~\cite{LPS88},
the  lower bounds on the girth of
regular graphs of degree $d\ge 10$ not a prime power
were largely increased.
Indeed, is obtained $\gamma_d \ge 1,06$
and even $\gamma_d \ge 1,33$ for larger values of $d$.
This improves upon the 30 years old $\gamma_d \ge 0,48$ proved
in~\cite{Im84}, for $d\not=2^k+1$. For $d=2^k+1$, this improves upon the $\gamma_d\ge \frac 2 3$
of~\cite{Mo94}.
It even outperforms
what the probabilistic method~\cite{ErSa63}
is able to give, namely $\gamma_d\ge 1$.

The construction of Imrich~\cite{Im84} is inspired by the previous
work of Margulis~\cite{Ma82}.
The families
that are built therein are derived from a motherh graph, seen as a Cayley graph on a suitable
free subgroup of $SL_2(\Z)$.
This prevents to use quaternions as done here and in~\cite{Ma88,Ma82,LPS88},
because the Hamilton quaternion algebra $\H(\Q)$ is not split (no isomorphism with the 2-by-2
matrices). Thanks to quaternions, it is comparatively possible to do better. The lower bound
obtained on the girth of the non-bipartite Cayley graphs $H_n$ on $PSL_2(\F_q)$  in Theorem 1.1, is
$\ge \frac{1,33}2 
\cdot \log_d |H_n|$ for $d$ large enough. As already mentioned, this is better than for the Cayley
graphs on $SL_2(\F_q)$ in~\cite{Im84}, where the lower bound on the girth is worked out directly on matrices
of $SL_2(\Z)$ (see Proposition 4 of~\cite{BoGa08} for more details) and not on integral quaternions as done
here.

\paragraph{Expander graphs}
It should be mentioned that all the families
of non-bipartite $d+1$-regular graphs $\Y_{d,p}$
defined in~\eqref{eq:Y}
are {\em expander families}.
This is due to their large girth property,
for which the theorem of Bourgain \& Gamburd~\cite[Theorem~3]{BoGa08}
holds.
In particular, the non-bipartite graphs $G_{d,p,q}$ do not have
a small chromatic number,  but have a small diameter in the order
of $O(\log |G_{d,p,q}|)$ (see~\cite[pp.455]{HoLiWi06}).

\paragraph{About possible generalizations}

In 1994, Morgenstern in~\cite{Mo94} 
has extended the construction of families
of $p + 1$-regular Ramanujan graphs by Lubotzky-Philips-Sarnak~\cite{LPS88} and Margulis~\cite{Ma88}
coming with a construction of families of $p^k + 1$-regular graphs, $p$ any prime and $k \in \N^\star$. The
idea was to use quaternion algebras over function fields that are of class number equal to 1
(admit a unique factorization property similar to Theorem 2.2).
Applying the technique developed in the present paper
to those graphs raises the hope to improve further more
 the estimates on the girth: 
Indeed, given an integer $d$ the next prime power $p^k$
 is always smaller than the next prime $p'$: $p^k \le p'$ (remember
that this ``gap'' plays an important role in the estimate of the girth). However,
Dickson's result do not hold directly
for the group $PSL_2(\F_{p^k} )$ and thus cannot guarantee the connectedness 
of the graphs as was done here. 
We have not tried, but even if connectedness
can be obtained in some cases,
we found out that the use of Morgenstern graphs
may not be worth
considering the tradeoff between simplicity and sharpness of the bounds, as explained below: 
\begin{itemize}

\item for an even number $d$, to build a $d + 1$-regular tree was required some ``involutions'' in
$\Pp(p)$, as explained in Remark~\ref{rem:free-group}. They were proved to exist only if $p\equiv 3 \bmod 8$. There
is no such involution in the similar special set of prime quaternions of Equality (9) of~\cite{Mo94}
(see Definitions 4.3 and 4.6 therein). Hence, to build a $d + 1$-regular tree we are led to
consider the prime $p = 2$, and to choose $d+1$ elements in the set defined in Equality (18)
and Definition 5.3 of~\cite{Mo94} (indeed, by Corollary 5.7 they yield such involutions). But in
this case, roughly because $PSL_2(\F_{2^k} ) = PGL_2(\F_{2^k} )$, the Cayley graphs $\Gamma_g$ obtained are
non-bipartite and only of girth $\ge \frac 2 3 \log_q |\Gamma_g|$ 
(see Theorem 5.13). This does not compete
with the girth of the graphs described in the present paper, even in the non-bipartite case.
\item for an odd number $d$, 
the use of Morgenstern graphs could make sense if connectedness is proved, however the values
of $c(d)$ 
for $d$ odd shown in Theorem~\ref{th:contrib} are not too bad, becoming close to the upper limit
$\frac 43$  rather quickly.

\item the use of the construction of Morgenstern would induce a jump in technicality, with additional new
results to address the problem of connectedness, and
without
a significant strengthening of the results, as shown by the two previous points.
\end{itemize}

\bigskip
\centerline{\bf Acknowledgment}

I am indebted  to J.-P. Tillich who initiated me to the subject
of Ramanujan graphs.
I would also like to thank the anonymous referees for their careful proofreading
and valuable comments.


\begin{thebibliography}{10}

\bibitem{Bo78}
B.~Bollob{\'a}s.
\newblock {\em Extremal graph theory}, volume~11 of {\em London Mathematical
  Society Monographs}.
\newblock Academic Press Inc. [Harcourt Brace Jovanovich Publishers], London,
  1978.

\bibitem{BoGa08}
J.~Bourgain and A.~Gamburd.
\newblock Uniform expansion bounds for {C}ayley graphs of {$SL_2(\F_p)$}.
\newblock {\em Ann. of Maths}, 167(2):625--642, 2008.

\bibitem{DaSaVa03}
G.~Davidoff, P.~Sarnak, and A.~Valette.
\newblock {\em Elementary number theory, group theory, and {R}amanujan graphs},
  volume~55 of {\em London Math. Soc. Student Texts}.
\newblock Cambridge U. Press, 2003.

\bibitem{Du99}
P.~Dusart.
\newblock The {$k$}-th prime is greater than {$k(ln k + ln ln k -1)$} for
  {$k\ge 2$}.
\newblock {\em Mathematics of Computation of the American Mathematical
  Society}, 68(225):411--415, 1999.

\bibitem{ErSa63}
P.~Erd{\"o}s and H.~Sachs.
\newblock Regul{\"a}re {G}raphen gegebener {T}ailenweite mit minimaler
  {K}nollenzahh.
\newblock {\em Wiss. Z. Univ. Halle-Willenberg Math. Nat.}, 12:251--258, 1963.

\bibitem{GaHoShShVi09}
A.~Gamburd, S.~Hoory, M.~Shahshahani, A.~Shalev, and B.~Vir\'ag.
\newblock On the girth of random {C}ayley graphs.
\newblock {\em Random Structures and Algorithms}, 35(1):100 -- 117, 2009.

\bibitem{Gr95}
A.~Granville.
\newblock Harald {C}ram{\'e}r and the distribution of prime numbers.
\newblock {\em Scandinavian Actuarial Journal}, 1995(1):12--28, 1995.

\bibitem{HoLiWi06}
S.~Hoory, N.~Linial, and A.~Wigderson.
\newblock Expander graphs and their applications.
\newblock {\em Bull. Amer. Math. Soc. (N.S.)}, 43(4):439--561 (electronic),
  2006.

\bibitem{Im84}
W.~Imrich.
\newblock Explicit construction of regular graphs without small cycles.
\newblock {\em Combinatorica}, 4(1):53--59, 1984.

\bibitem{LaUs95}
F.~Lazebnik and V.~A. Ustimenko.
\newblock Explicit construction of graphs with an arbitrary large girth and of
  large size.
\newblock {\em Discrete Appl. Math.}, 60(1-3):275--284, 1995.
\newblock ARIDAM VI and VII (New Brunswick, NJ, 1991/1992).

\bibitem{LPS88}
A.~Lubotzky, R.~Phillips, and P.~Sarnak.
\newblock Ramanujan graphs.
\newblock {\em Combinatorica}, 8(3):261--277, 1988.

\bibitem{Ma82}
G.~A. Margulis.
\newblock Explicit constructions of graphs without short cycles and low density
  codes.
\newblock {\em Combinatorica}, 2(1):71--78, 1982.

\bibitem{Ma88}
G.~A. Margulis.
\newblock Explicit group-theoretic constructions of combinatorial schemes and
  their applications in the construction of expanders and concentrators.
\newblock {\em Problemy Peredachi Informatsii}, 24(1):51--60, 1988.

\bibitem{Mo94}
M.~Morgenstern.
\newblock Existence and explicit constructions of {$q+1$}-regular {R}amanujan
  graphs for every prime power {$q$}.
\newblock {\em J. Combin. Theory Ser. B}, 62(1):44--62, 1994.

\bibitem{RaRu96}
O.~Ramar\'e and R.~Rumely.
\newblock Primes in arithmetic progressions.
\newblock {\em Mathematics of Computation}, 65(213):397--425, 1996.

\bibitem{Se77arithm}
J.-P. Serre.
\newblock {\em Cours d'arithm{\'e}tique}, volume~2.
\newblock Presses universitaires de France, 1977.

\bibitem{Ta81}
R.~M. Tanner.
\newblock A recursive approach to low complexity codes.
\newblock {\em IEEE Trans. on Inform. Theory}, 27(5):533--547, 1981.

\end{thebibliography}
\end{document}